\documentclass[12pt,leqno]{amsart}
\usepackage{amsaddr}
\usepackage{geometry}                
\usepackage[utf8]{inputenc}
\usepackage{amstext,amsmath,amssymb,amsbsy,bm,amsthm}
\usepackage{color}
\usepackage{graphicx}
\def\Rd{\color{red}} 
\def\B{\color{black}}
 
\def\Bl{\color{blue}}

\numberwithin{equation}{section}

\newtheorem{definition}{Definition}[section]

\newtheorem{problem}[definition]{Problem}
\newtheorem{conjecture}{Conjecture}

\DeclareMathOperator*{\argmin}{arg\,min}

\usepackage{mathrsfs}

\title{Deep Learning for the Approximation of a Shape Functional}
\author[F. Calabrò  \and S. Cuomo \and F. Giampaolo \and S. Izzo \and alt.]{Francesco Calabrò \and Salvatore Cuomo \and Fabio Giampaolo \and Stefano Izzo \and Carlo Nitsch \and Francesco Piccialli \and Cristina Trombetti}
\address{Dipartimento di Matematica e Applicazioni ``R.Caccioppoli"\\ Università  degli Studi di Napoli Federico II}
\date{}
\begin{document}

\maketitle
\begin{abstract}
Artificial Neuronal Networks are models widely used for many scientific tasks. One  of the well-known field of application is the approximation of high-dimensional problems via Deep Learning. In the present paper we investigate the Deep Learning techniques applied to Shape Functionals, and we start from the so--called Torsional Rigidity. Our aim is to feed the Neuronal Network with digital approximations of the planar domains where the Torsion problem (a partial differential equation problem) is defined, and look for a prediction of the value of Torsion. Dealing with images, our choice fell on Convolutional Neural Network (CNN), and we train such a network using reference solutions obtained via Finite Element Method. Then, we tested the network against some well-known properties involving the Torsion as well as an old standing conjecture. In all cases, good approximation properties and accuracies occurred.
\end{abstract}


\section{Introduction}
By ``Shape Functional" one usually denotes any abstract functional defined on some class of subsets of $\mathbb{R}^d$. The Volume and the Perimeter are probably the most common examples of shape functionals. More sophisticated examples coming from physics include: drag force of a solid moving into a fluid, the electrostatic capacity of a conductor body in vacuum space, the effectiveness of the thermal insulation around a thermally conductor body etc. Very often (as in the examples above), there is some region (say $\Omega$) to which we associate a certain scalar value (say $F(\Omega)$).

Shape optimization is nowadays a very active field of research, ranging from abstract mathematics to applied physics and engineers. In a very broad sense, it aims at finding sets (called optimal sets) that minimize or maximize some shape functional, possibly under feasible geometrical constraints.

In many cases, in order to determine the value of a functional shape, one has to solve a differential problem: partial differential equation (PDE) or systems of PDEs.

Depending on the problem, the determination of the numerical solution can sometimes be very time-consuming. In addition, the class of domains in which one looks for the optimal sets is at best (unless one introduces some severe restrictions) infinite-dimensional, and there is no general method allowing us to find the optimal domain in a straightforward and reliable way.

In this paper, we investigate the possibility to train a Neural Network (NN) to approximate a shape functional. It can be considered as a first step toward the ultimate goal of a general implementation of NNs in shape optimization.

We warn the reader that the use of NNs is not meant to replace numerical methods, such as finite elments or finite difference schemes. Classical methods are still by far more accurate. They are actually used to train the network, since in absence of an analytical computation, they provide the only reference values. Neural Network can be used to support, to reinforce or to validate a numerical investigation. At the same time NNs are extremely efficient in terms of computational time. Some of the problems arising in shape optimization, expecially those involving a free boundary problem, are so complicated that classical numerical schemes are sometimes very difficult to implement and very time consuming. There is a good chance that a Deep Learning approach could be succesful in providing quick answers, which are very welcome, even at the risk of being suboptimal.


As a first attempt, we have chosen a functional called Torsion (or Torsional Rigidity). As we will see later (and as the name suggests), it is related to the mechanical response of an elastic body under Torsion. To determine its value, one has to solve a Poisson problem in a planar set $\Omega$. The geometry of $\Omega$ uniquely determines the value of the functional. Torsional rigidity has the great advantage of having been intensively investigated in the last decades. Many theoretical results, as well as some long standing conjecture, provide reliable tests to check the accuracy of the neural network. However, no matter how large is the training set, there is no way it could be a good representative of the large class of open subset of the plane. Our ultimate test was to check the ability of the NN to extrapolate information on {\it unseen} cases, those with topological properties significantly different from the set used during the training.\\


For the reader not accustomed to the NNs, the naive idea beyond our approach can be easily summarized as follows. Take any planar bounded set enclosed in a box. Such a set can be represented by an digital b/w image: a square matrix of zeros and ones. The value $1$ corresponds conventionally to a point inside the planar set and the value $0$ corresponds to a point outside the set. We aim at an unknown function, {\it Target Function}, (which certainly exists even though it may be extremely complicated) that takes all these ones and zeros and provides the real value corresponding to the torsional rigidity of our input set. We want to approximate the target function and a NN is indeed a good candidate to do this job. How?\\

The building brick of a neural networks is the ``neuron". A neuron itself is a function. It takes as input real values $x_j$ ($j=1,2,...,K$) 
and makes a linear combination $z$ defined as follows
\begin{equation*}
    z = \sum_{j}^K w_j x_j +b.
\end{equation*}
Here $w_j$ are called {\it weights}, and $b$ is called {\it bias}.
The output is then passed to a nonlinear activation function $g$. When 
\begin{equation*}
    g(z) = \dfrac{1}{1+e^{-z}}, 
\end{equation*}
we have the so--called {\it sigmoid neuron}.
Weights and bias have to be fixed, and they are called {\it trainable parameters} of the neuron.
A different activation function provides a different type of neuron. For instance the well--known {\it RELU} (Rectified Linear Unit) has the following activation function
\begin{equation*}
    g(z) = \max\{0,z\}. 
\end{equation*}
All neurons share the common feature of taking a certain number of inputs and providing an output. 
However, a single Neuron is ``just'' a simple function; in order to extract complex information from the input, it is necessary to combine multiple neurons, in such a way that the output of a neuron can become the input of other neurons and so on. Each weight then represent the strenght of the connection between two neurons, reminiscent of our brain neuronal structure.\\
What in general is done is to organize the neurons in structures, called \textit{layers}, to augment the ability to extract patterns of data and then concatenate layers to elaborate these patterns. 
Such an architecture, composed of many connected neurons, when fully connected, is called Deep Neural Network (DNN).
At the end, the output of this structure is a composition of neurons' functions, whose distance to the ``target function" is measured by a ``loss function". The structure which we consider is classified as feedforward neural network. The layers are ordered and the information flows from the input layer, to the output layer. No loop is formed between the neurons.

In a NN, weights and bias have to be changed to minimize the loss function. Such an operation is the most delicate, but there are some standard procedures for the optimization, including a particular gradient descent procedure called backpropagation \cite{lecunn1988atheoretical}. This learning procedure provides an automatic strategy for updating the trainable parameters (weights and bias) to obtain a good fitting model for the considered input data. \\
The interested reader should look for Universal Approximation Theorems \cite{cybenko1989approximation} which justify the use of NN to ``represent'' a wide variety of unknown nonlinear functions.


Over the years several types of NNs have become very popular. Among the others, important examples are the Convolutional Neural Networks (CNN), which are widely used to extract spatial information from data taking advantage of the property of space invariance of the Convolution operation. In particular, CNNs are mostly desiged to deal with images. 

We shall see that a trained CNN can approximate very well our dataset. In particular, we gain good approximation error - thus accuracy - with a computational cost of evaluating Torsion that is orders of magnitude lower when compared to the standard resolution. The cost of the overall approximation is moved to the training of the CNN that is done only once and in advance. The resulting CNN is able to compute the Torsion in real-time. 
\\
The paper is organized as follows. In Section \ref{sect:Problem} we present the formal mathematical definition of Torsional Rigidity together with some known theoretical results. In Section \ref{sect:Model} we describe the construction of the training data set. We first present the routine used for the random generation of planar domains and then, the numerical approach (Galerkin method) for solving the PDE problem on such domains. In Section \ref{sect:Deep} we describe the Deep Neural Network that mimes the torsional rigidity. The accuracy of the DNN on training, validation and test set is discussed in Section \ref{sect:Test}. We then focus on polygonal domains  
and we show that, with a fairly good accuracy, the NN ``agrees'' with a long standing conjecture. 
Finally, in Section \ref{sect:Predictions} we test the NN on some complicated examples which are topologically different from the training images. First we ask the NN to guess the Torsion on disconnected domains. Then we try out the tricky case of multiple connected sets. In both cases the results are very satisfactory. The paper ends (Section \ref{conclusioni}) with some conclusions and perspectives.

\section{The considered problem}\label{sect:Problem}
In Mathematical Analysis the Torsional Rigidity $T(\Omega)$ of a bounded open set $\Omega\subset \mathbb{R}^d$ is defined as
\begin{equation}\label{variazionale}
\displaystyle T(\Omega)^{-1}=\min_{v\in H^1_0(\Omega)} \frac{\displaystyle\int_\Omega|\nabla v|^2 \,dx}{\left(\displaystyle\int_\Omega v \,dx\right)^2} 
\end{equation}

It is named after the physical interpretation, in case $\Omega$ is a planar set ($d=2$), of representing, in the framework of elasticity, the torsional rigidity of a prismatic bar of cross section $\Omega$. \\
A minimizer $u$ of \eqref{variazionale} is provided by the unique solution to the following problem
\begin{equation}\label{PDE}
\left\{\begin{array}{ll}
-\Delta u=1 &\mbox{in $\Omega$,}\\\\
u=0 &\mbox{on $\partial \Omega$.}
\end{array}\right.
\end{equation}
The function $u$ is the so-called Torsion Function (or Prandtl stress function, see \cite{TimoshenkoStephen1951Toe}) and
the torsional rigidity becomes (using equations \eqref{PDE}--\eqref{variazionale} together with the divergence theorem)
$$\displaystyle T(\Omega)=\int_\Omega u \,dx.$$

The Torsional Rigidity (or just Torsion), is uniquely determined by the geometry of $\Omega$ and it is one of the most simple examples of so-called \emph{Shape Funcional}. The dependence of the functional from the set $\Omega$ is provided by the solution to a partial differential equation and therefore non trivially deductible from the geometrical properties of $\Omega$.
The Torsion has been intensively investigated and although explicit solutions are known only in very few cases, a large number of known properties can be found in the literature. We summarize here the most important features.
\begin{itemize}
\item (Rototranslation) The Torsion is invariant under rototranslation.
\item (Scaling) Under scaling the Torsion changes according to a power law: when the set $\Omega$ is replaced by an omothetic copy $t\,\Omega$ (for some $t>0$), then $T(t\Omega)=t^4\,T(\Omega)$.
\item (Additivity) Torsion is additive when acting on unions of disjoint open sets: $T(\Omega)=T(\Omega_1)+T(\Omega_2)$ whenever $\Omega_1$ and $\Omega_2$ are open, disjoint and such that $\Omega_1 \cup \Omega_2 = \Omega$.
\item (Monotonicity) As a trivial consequence of the maximum principle for Poisson problem, the Torsion is monotone by inclusion: $T(\Omega_1)\le T(\Omega_2)$ whenever $\Omega_1\subset\Omega_2$.
\item (Saint-Venant) The Torsion of a set $\Omega$ is always smaller then the one of a ball having same volume.
\end{itemize}

The first four properties follow trivially from the definition. The fifth one is named after Saint-Venant, who first conjectured it back in 1856. It was proved by P\'olya in 1948 (see \cite{PolyaSzegoBook1951,PolyaTorsion1948}). It means that among beams with constant cross sectional area, the  maximal Torsional Rigidity is achieved by circular cylinders.\\
Observe that, from property (Scaling), denoting by $A(\Omega)$ the area of an  open set $\Omega$, then $A^{-2}(\Omega)T(\Omega)$ turns out to be scaling invariant. Indeed, in the following sections, we will often plot the Torsion $T$ against $A^2$. 

Moreover, if $\mathcal{D}$ denotes any disk, then Property (Saint-Venant) reads 
\begin{equation}\label{eq:upperboundBall}
    A^{-2}(\Omega)T(\Omega)\le A^{-2}(\mathcal{D})T(\mathcal{D}).
\end{equation}

The elementary form of the solution $u$ to \eqref{PDE} is known only for very few domains. Amog the others:
\begin{enumerate}
    \item The disk $\mathcal{D}_R$ of radius $R$ where
    \begin{equation*}
        u(x,y)=\frac14 (R^2-x^2-y^2)
    \end{equation*}
    and
    \begin{equation}\label{eq:TorBall}
        T(\mathcal{D}_R)=\frac{\pi}8\, R^4 
    \end{equation}
    \item More generally, on Ellipse $\mathcal{E}_{ab}$ of equation $\frac{x^2}{a^2}+\frac{y^2}{b^2}-1<0$ where
    \begin{equation*}
        u(x,y)=\frac{a^2b^2}{2(a^2+b^2)} \left(1-\frac{x^2}{a^2}-\frac{y^2}{b^2}\right)
    \end{equation*}
    and
    \begin{equation}\label{eq:TorEllipse}
        T(\mathcal{E}_{ab})= \frac{2a^3b^3}{a^2+b^2}\, \frac{\pi}{8}
    \end{equation}
     \item The annular ring $\mathcal{A}_{rR}$ of equation $r^2<x^2+y^2<R^2$ where
    \begin{equation*}
         u(x,y)= -\dfrac{x^2+y^2}{4} + A \log (x^2+y^2) + B 
    \end{equation*}
    (the constants $A$ and $B$ are determinated by the conditions $u(R)=u(r)=0$) and
    \begin{equation}\label{eq:TorAnnularRing}
        T(\mathcal{A}_{rR})=\frac{\pi } {8}\, \,  \left( (R^4-r^4) - \frac{\left( R^2 - r^2 \right)^2}{\log(R/r) } \right)
    \end{equation}
\end{enumerate}
For equilateral triangles it is possible to write down an explicit solution but already for rectangle, the solution can at best be expressed as a series of elementary functions and the value of the Torsion as well requires the summation of a series \cite{TimoshenkoStephen1951Toe}. Already for general polyogns the solutions in closed elementary form are out of question. 
A famous conjecture reads as follows
\begin{conjecture}\label{congettura}
Among all pentagons of given area the regular pentagon achieves the maximal Torsional rigidity 
\end{conjecture}
Actually the original conjecture has been stated for the Laplace eigenvalue (see forn instance \cite{HenrotBook2006}), but it is generally tacitally assumed to refer also to the Torsional rigidity.
Similar conjectures stand still unsolved for any number of sides (greater than $5$). For triangles and quadrilaterals it is easy to prove that the property holds true (\cite{HenrotBook2006}).



\section{Model the Torsion approximation}\label{sect:Model}

\begin{figure}
	\centering
	\includegraphics[width=\textwidth]{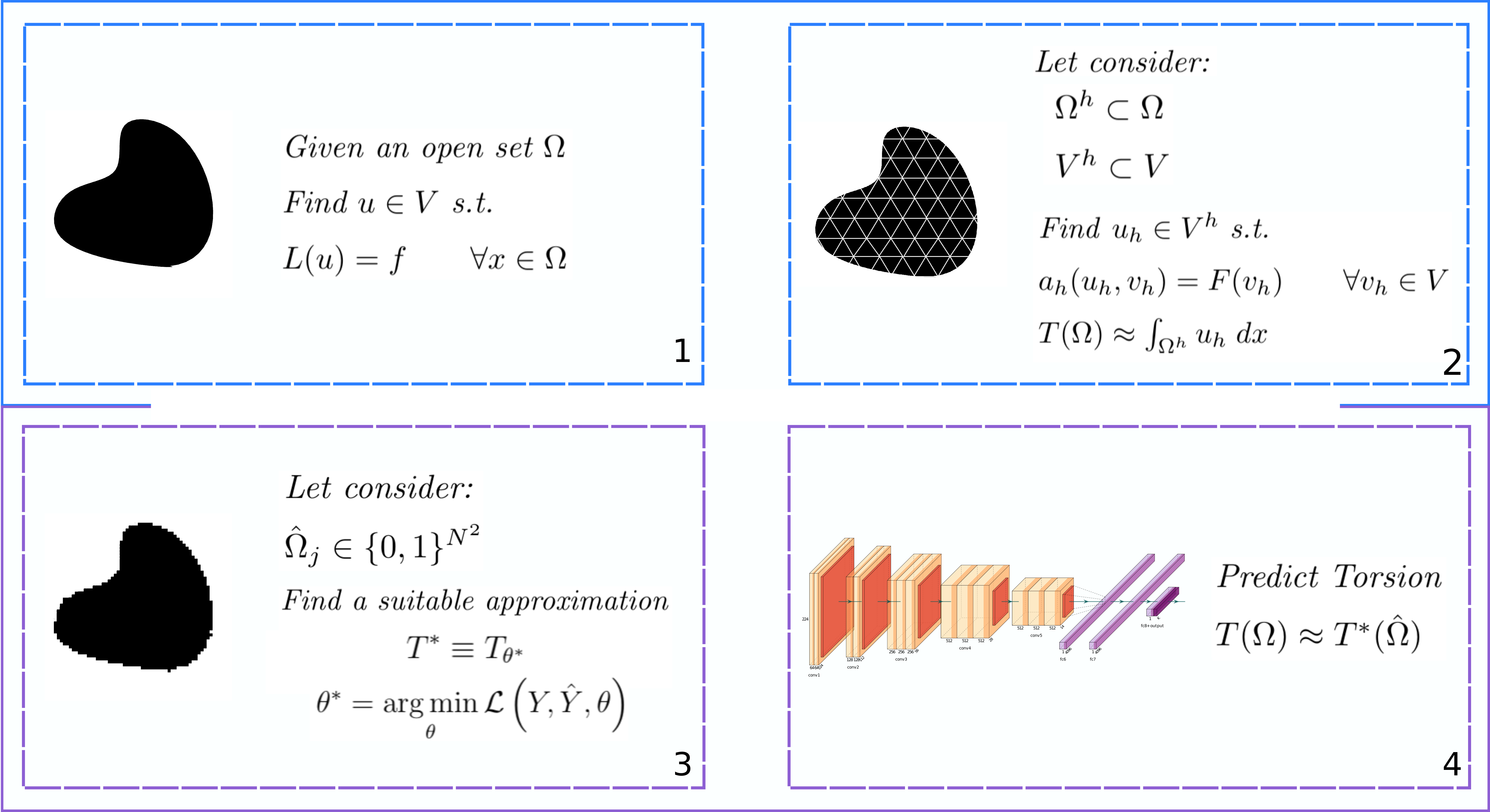}
	\caption{ The proposed DL pipeline is composed of four stages. Starting from the open set $\Omega \subset \mathbb{R}^2$ and the matematical problem  (subfigure \textbf{1}), such problem is numerically solved by means of a Galerkin method (subfigure \textbf{2}). These procedure is repeated for different sets $\Omega$ in order to populate a reference database. Conversely, a DL approximation model to compute the Torsion is identified (subfigure \textbf{3}) by adopting a supervised approach from a binary image $\hat{\Omega}$ of $\Omega$. Through the DL model, we obtain the approximated Torsion value of $\hat{\Omega}$ (subfigure \textbf{4})}
	\label{fig:pipeline}
\end{figure}

In this section we describe the model used to predict $T(\Omega)$ by means of a trained DNN. 
More specifically, in Subsection \ref{subsect:Galerkin} we discuss how to construct a training dataset and in Subsection \ref{subsect:approximation} we describe the ML methodology employed to approximate the functional $T(\cdot)$.

\subsection{Training dataset via Galerkin method}\label{subsect:Galerkin}
The first step was the construction of different instances of $\Omega$. More specifically $\{\Omega_j, \, j = 1,2,...,M\}$ will denote a set of $M$ open bounded subset of $\mathbb{R}^2$ randomly generated. For simplicity, we decided to confine all sets into a reference square $\mathcal{Q}=[-2,2]^2$. In order to generate a domain $\Omega$ we selected a random number $r$ (between $3$ and $20$) and a random set of $r$ points into $\mathcal{Q}$. Starting from such a set of points, $\Omega$ was constructed by means of interpolation using \texttt{polyshape} routine in MATLAB and, in some randomly selected cases, the \texttt{spline} routine (to smooth the boundary). Therefore the boundary of $\Omega$ was either a closed polygonal path that passes through (a subset of) the random points, or a closed curve interpolating (a subset of) the random points with a spline. \\ 
Thereafter we considered problem \eqref{PDE} on each of the previous set and we looked for a numerical approximate value of the Torsion. 
In a fairly standard way, the domain $\Omega$ was approximated by a computational domain $\Omega^h$. On such a domain, a weak formulation of the original PDE was considered. By the Galerkin method \cite{quarteroni2009numerical} we found an approximate solution $u_h$ to the original problem \eqref{PDE} written as linear combination of basis functions $\phi_i$: 
\[ u^h\in V^h:=span\{\phi_i\} .  \] 
Here $V^h$ is the space of piece--wise linear functions, leading to the so-called P1 Finite Elements (FE).
We repeat the computations on $M$ domains $\Omega_j, j=1,\dots, M$ and denote the approximate value of the Torsion  by $\bar{T}({\Omega_j})$. 
\\
All the computations were made using Matlab, in particular the \texttt{PDEs toolbox} for the FE, the core \texttt{Computational Geometry} repository for the construction of the different istances of $\Omega_j$, and \texttt{solvepde} routine to solve numerically the PDE. \\

By $Y= \{y_j=\bar{T}({\Omega_j}) \in \mathbb{R}^+ \text{ with }  j=1,...,M \}$ we shall denote the set of Torsion values obtained solving numerically the PDE problem \eqref{PDE}.\\

\begin{figure}\label{fig:domini}
	\centering
	\includegraphics[width=\textwidth]{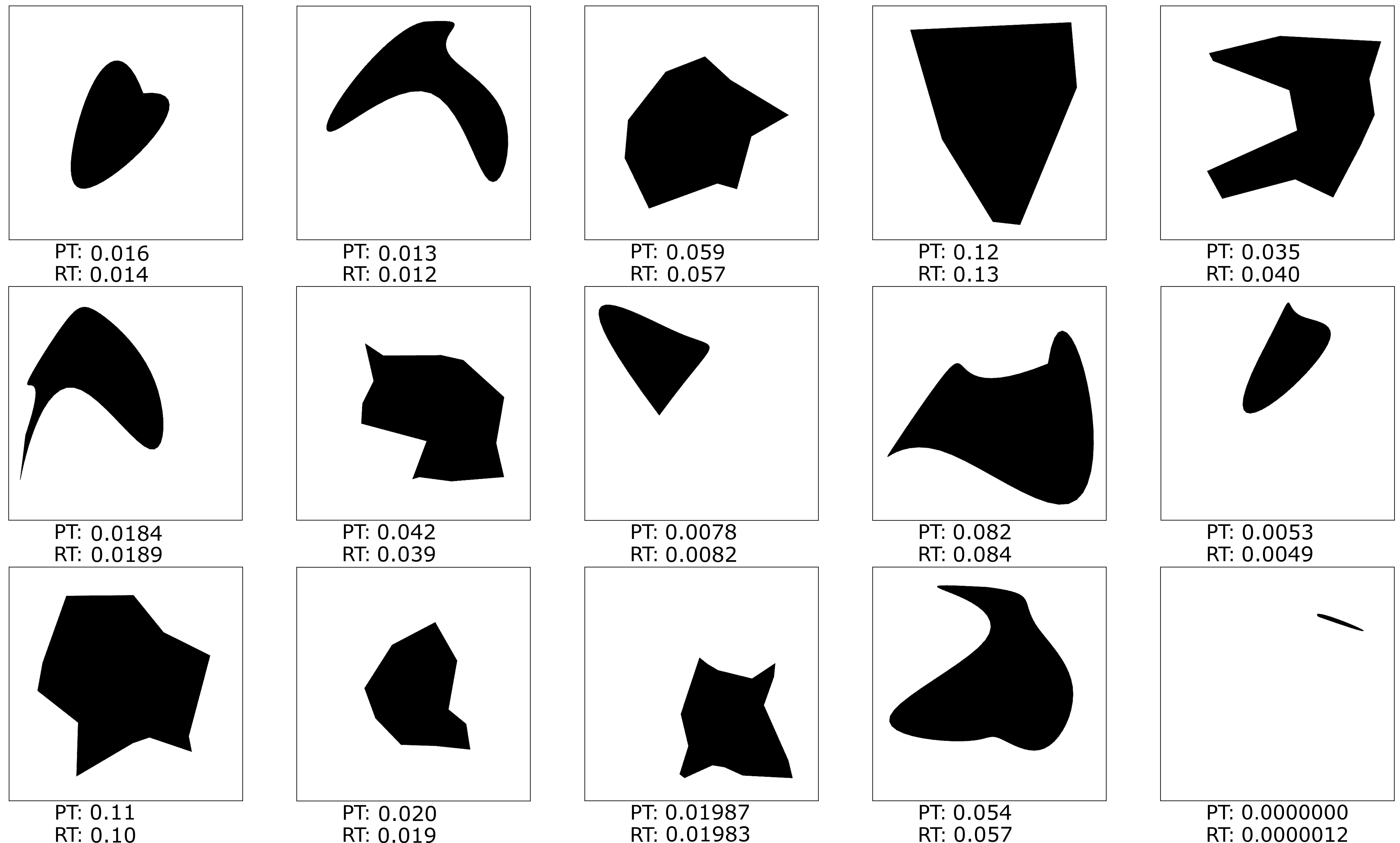}
	\caption{In this figure, some of the constructed domains used to populate the traing set. With PT we denote the value predicted by the DNN. With RT the reference solution for Torsion, numerically computed.}
\end{figure}

In Figure \ref{fig:domini}, we present some of the domains generated and the numerical value of the Torsion (RT {\it Reference Torsion value}) . \\

\subsection{Computation of Torsion via NN}\label{subsect:approximation}
The first step consisted in turning $\{\Omega_j, \, j = 1,2,...,M\}$ into a reasonable binary input. A straightforward way is to use an image for each domain $\Omega_j$. Black and white images of size $N\times N$ pixels, basically a square matrix of zeros and ones. \\
By $\hat{\Omega}_j \in \{0,1\}^{N^2}$ with $j = 1,2,...,M$ we denote the binary images corresponding to the sets $\Omega_j$.\\

Then we set up the approximation problem as follows:
\begin{problem}
	Given a finite \textbf{set of instances} $\hat{\Omega}=\{\hat{\Omega}_j\}_{j\in \{1,...,M\}}$ and an M-tuple $y\in\mathbb{R^+}^M$ called \textbf{target}, the set of M pairs (\textbf{samples}) $\{(\hat{\Omega}_j, y_j)\}_{j\in\{1,...,M\}}$ is defined \textbf{training set}. \\
	\textbf{The regression problem} consists in finding a function (model)
	\[ \hat{T}:\hat{\Omega}_j\in \{0,1\}^{N^2} \to \hat{T}(\hat{\Omega}_j)\in \mathbf{R}^+ \]
	such that
	\[ \hat{T}(\hat{\Omega}_j)=\bar{T}({\Omega}_j) \ . \]
\end{problem} \qed \\

The first question which arises is how to choose the space of approximating functions. The  problem is a high dimensional approximation problem: in view of the classical Weierstrass Theorem one can think to use standard procedure to construct such approximation. Indeed, such problem cannot be solved pratically using classical approaches such as polynomials, piecewise polynomials, or truncated Fourier series because of the so called curse of dimensionality: 
the number $N_{dof}$ of parameters required to achieve certain accuracy depends exponentially on the dimension $N^2$.
\\
On the other hand, the Universal Approximation
Theorem states that under very mild conditions, any continuous function can be uniformly approximated on compact domains by neural network functions, see \cite{hornik1990universal,kratsios2019universal,mallat2016understanding,mhaskar2016deep,petersen2018optimal,zhou2020universality}. In the recent work \cite{zhou2020universality} the course of dimensionality has been overcome in the case of deep CNN. These are the common networks used for the treatment of problems involving images \cite{jin2017deep,unser2019representer}
for instance in medical applications such as in automatic tumor markers \cite{piccialli2021survey,plis2014deep}.\\
For all these considerations, we decided to use a set of deep CNN (see next section for more details).

Our numerical problem then is to find a suitable CNN 
that depends on a finite set of (trainable) parameters $\theta\in\mathbf{R}^{N_{dof}}$. Since $\hat T(\cdot)$ is uniquely determined by the training set of parameters $\theta$, when necessary, we explictly emphasize such a dependence by writing $\hat T_\theta(\cdot)$. \\
The optimization reads:
\begin{equation}\label{eq:approx_tors_funzione}
\text{Find }    \theta^*\in\mathbf{R}^{N_{dof}} \text{ such that } {T}^*\equiv\hat{T}_{\theta^*}
\approx \bar{T}. 
\end{equation}
The condition to be satisfied is that
$T^*(\hat{\Omega}_j)$ is close to $\bar T({\Omega}_j)$ $\forall j\in\{1,...,M\}$ according to a certain cost function $\mathcal{L}$ (called also loss function). Namely, if we denote by $\hat{Y}= \{\hat{y}_j=\hat{T}({\hat{\Omega}_j}) \in \mathbb{R}^+ \text{ with }  j=1,...,M \}$ the set of Torsion values obtained by the neural network, we want to minimize a suitable notion of {\it distance} between $Y$ and $\hat Y$.\\

We generated $5\cdot 10^3$ 
images of size $677\times 677$ pixels and calculate reference Torsion values on each of these: the target values. To such a set we applied a simple procedure of data augmentation which consists in considering new copies obtained by rotating and translating the original images. We denote the resulting set, made of $2\cdot 10^4$ 
images, by $\mathcal{F}$. Then, we randomly splitted this set in three disjoint subsets: training set $\mathrm{Tr}$, validation set $\mathrm{Val}$ and test set $\mathrm{Te}$: 
\begin{equation} \label{def:testset}
   \mathcal{F} = \mathrm{Tr}\cup \mathrm{Val} \cup \mathrm{Te} \ .
\end{equation}
Following a standard procedure, $ \mathrm{Tr}, \; \mathrm{Val}, \; \mathrm{Te} $ are selected containing more or less 70\%, 10\% and 20\% of the images, respectively. Common practice is indeed to tune the parameter of the Network to match the training set. Use the Validation set for instance to take under control, during the optimization process, overfitting issues. At the end of the process, the test set provides an unbaised evaluation of the final NN, and does not influence directly the process of learning.

From now on we refer to $Y$ and $\hat Y$ to denote the training data and the corresponding values obtained by the neural network. Therefore $M$ represents just the number of images in the training set.\\
The optimal vaules $\theta^*$ are determinated by an optimization procedure that minimizes the loss function:
\begin{equation}\label{mse regolarizzato}
     \mathcal{L}(\hat{Y}, Y, \theta)=\dfrac{\sum_{i=1}^{M} \left(y_i-\hat{y}_i\right)^2}{M}+\lambda\sum_{j=1}^v w_j^2
\end{equation}
where 
$\{w_j\}_{j=1,...,v} $ are the weights of the network: a subset of the parameters $\theta$. Here $\lambda$ is some positive constant to be choosen in advance.\\

The first term on the righthand side of \eqref{mse regolarizzato} is the Mean Squared Error (MSE) between target and predicted values. The second term consists of the penalty term: a fixed positive parameter $\lambda$ multiplied by the squared $\ell^2$ norm of $\{w_j\}$. 
The MSE is an widely used loss function in Machine Learning problem; it exhibits robustness to outliers. The main drawback of this measure is that it loses accuracy on small data. 
On the other hand, the actual value of the Torsion is small each time the set has a very small area (remember that by omothety the Torsion rescales as the square of the area) and here is when the approximation by pixels becomes less and less accurate. Therefore there is no need to push the Network to look for accuracy in these cases.\\
Finally, the penalty (or regularization) term is a standard practice in ML as well. It allows us to control the amplitude of the network weights, and reduces the danger of an overfitting phenomena.\\

So, in summary, we want to determine

\begin{equation}\label{ML equation}
    \theta^* = \argmin_{\substack{\theta}} \mathcal{L}\left( Y, \hat Y, \theta \right)
\end{equation}
and
\begin{equation*}
    T^*\equiv T_{\theta^*}.
\end{equation*}

The minimization is obtained through a training process performed on the network called backward propagation method. Similarly to a gradient descent, the training process generates a sequence of DL neural networks that has to be stopped at the $k$-iteration. The number of iterations ({\it epoch}) $k$ has to be choosen accurately. Eventually, the parameters $\theta^{(k)}$ defines the network.\\


\section{Neural Network Architecture}\label{sect:Deep}

\begin{figure}
	\centering
	\includegraphics[width=\textwidth]{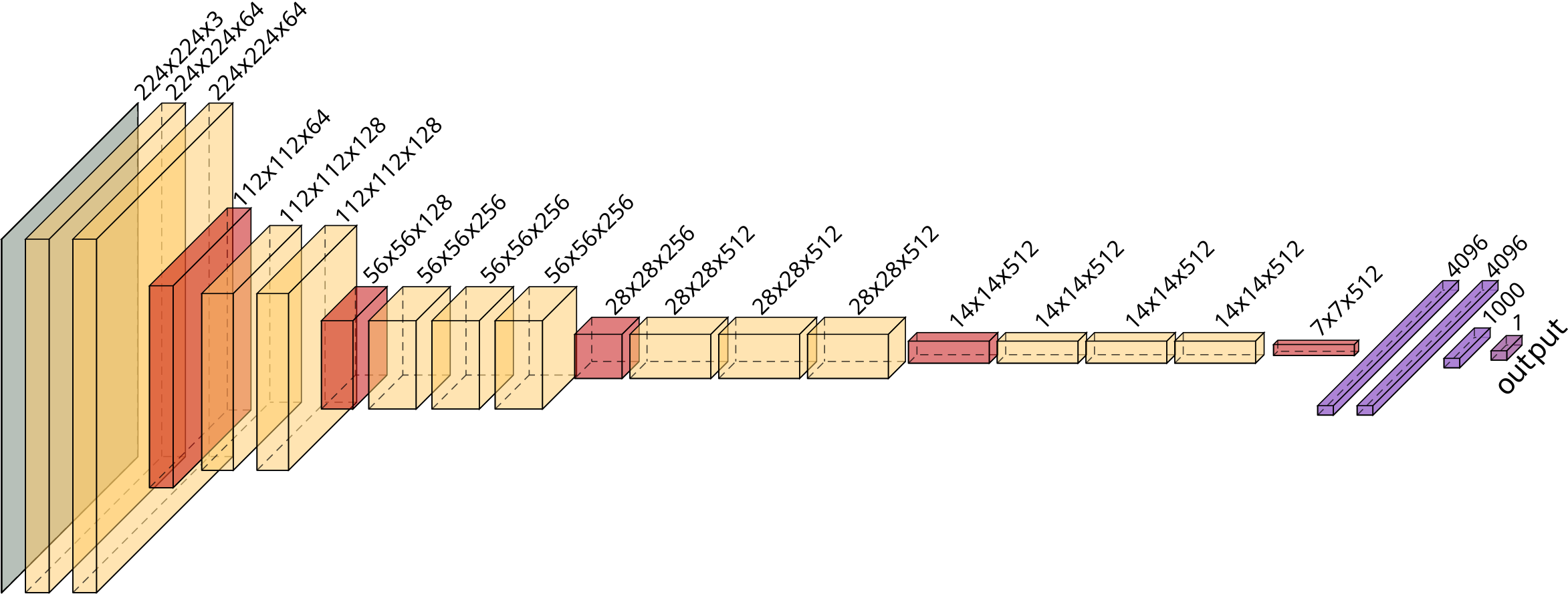}
	\caption{The VGG16-Net architecture. It is composed of 16 layers. The yellow layers are Convolutional layers; the red is the MaxPooling layer, the purple is the Fully Connected layer. In particular, the first five stacks are composed by $2$ or $3$ convolutional layers of $64, \, 128, \, 256, \, 512, \, 512$ kernels respectively with dimension $3\times 3$ respectively and a layer of a bidimensional MaxPooling at the end of each stack with pool size $2\times 2$. The second part of the network contains four Fully Connected layers with parameters $4096, \, 4096, \, 1000, \, 1$, respectively. The final layer is the output of the whole process.}
	\label{fig:VGG16}
\end{figure}

In this section, we describe the architecture of the CNN model (see Figure \textbf{\ref{fig:VGG16}}). The Convolutional Neural Network model takes as input an image $\hat{\Omega}$ and gives as ouput $T^*(\hat{\Omega})\approx T(\Omega)$.
The model applies discrete convolutional operations between the matrix representing the image $\hat{\Omega}$ and another matrix called kernel $K$. The result of this operation is called feature map $C$ 
\begin{equation}
    C(i,j) = (\hat{\Omega}*K) (i,j) = \sum_m \sum_n \hat{\Omega}(m,n)K(i-m,j-n) 
\end{equation}
What happens in a convolutional layer can be summarized in the following steps: random kernels are generated (with a fixed size), and convolutional operations are applied on all the inputs $\hat{\Omega}$ with the same $K$. 
A parameter called \textit{stride} manages the displacement of the filter along the image, i.e. how changes the index $(i,j)$. 
It is also possible to change the kernel and repeat the procedure, obtaining a different feature map. The set of feature maps obtained changing the kernel elements is the real output of the convolutional layer. \\
It is possible to concatenate the convolutional layer with other layers (also of different types), obtaining a multi-layer neural network called Deep Convolutional Neural Network (DCNN). In our case, we use a very well known DCNN called VGG16-Net (see Figure \ref{fig:VGG16}) \cite{simonyan2014very}. More in detail, the computational blocks highlighted in yellow are Convolutional layers; the red ones are MaxPooling kernels, i.e. a downsampling layer that extracts the max value from a fixed window. Moreover, the purple layer is a Fully Connected layer, i.e. a layer whose inside neurons connect to every neuron in the preceding layer. The input of VGG16-Net is an image with dimension $224\times 224$ pixel, so some steps of preprocessing consist of a procedure for image downscaling. The final output of the network is the predicted value of Torsion of the input image. 
Vgg16-Net is a Convolutional Neural Network composed of 13 convolutional layers and three fully connected layers with activation function ReLU \cite{lecun2015deep}. A peculiarity of this NN is the small kernel size, which is important to capture the basic spatial information. 
Furthermore, the zero-padding is used to preserve the dimensions. 
\B
\section{Experimental results on the trained network}\label{sect:Test}
After the training procedure, we tested the model in several ways.
\subsection{Test set: the network robustness}\label{subsec:test_set}
Following a stadard procedure, in the first step, we use the classical test set $\mathrm{Te}$ defined in \eqref{def:testset}, made of images generated with same method used for the training one. 
The test set contains more or less $4000$ images. 
The results are presented in Table \ref{tab:VGG-accuratezza} together with those for the training and validation. The parameter used are: learning rate of $10^{-4}$, regularization of $10^{-6}$, dropout of $0.5$ and batch size of $16$. \\
\B

\begin{table}[h]
    \centering
    \begin{tabular}{|l|c|c|} \cline{2-3} 
        \multicolumn{1}{c|}{} & Loss Function value & MSE \\ \hline 
        Tr (Training set) & $ 3\cdot 10^{-5} $ & $ 2\cdot10^{-5} $   \\ \hline
        Val (Validation set) & $1.6\cdot 10^{-5} $ & $8\cdot10^{-6} $  \\ \hline
        Te (Test set) & $1.4\cdot 10^{-5} $ & $6\cdot10^{-6} $  \\ \hline
    \end{tabular}
    \caption{Loss function and MSE of the VGG-16.}
    \label{tab:VGG-accuratezza}
\end{table}

\begin{figure}
	\centering
	\includegraphics[width=0.45\textwidth]{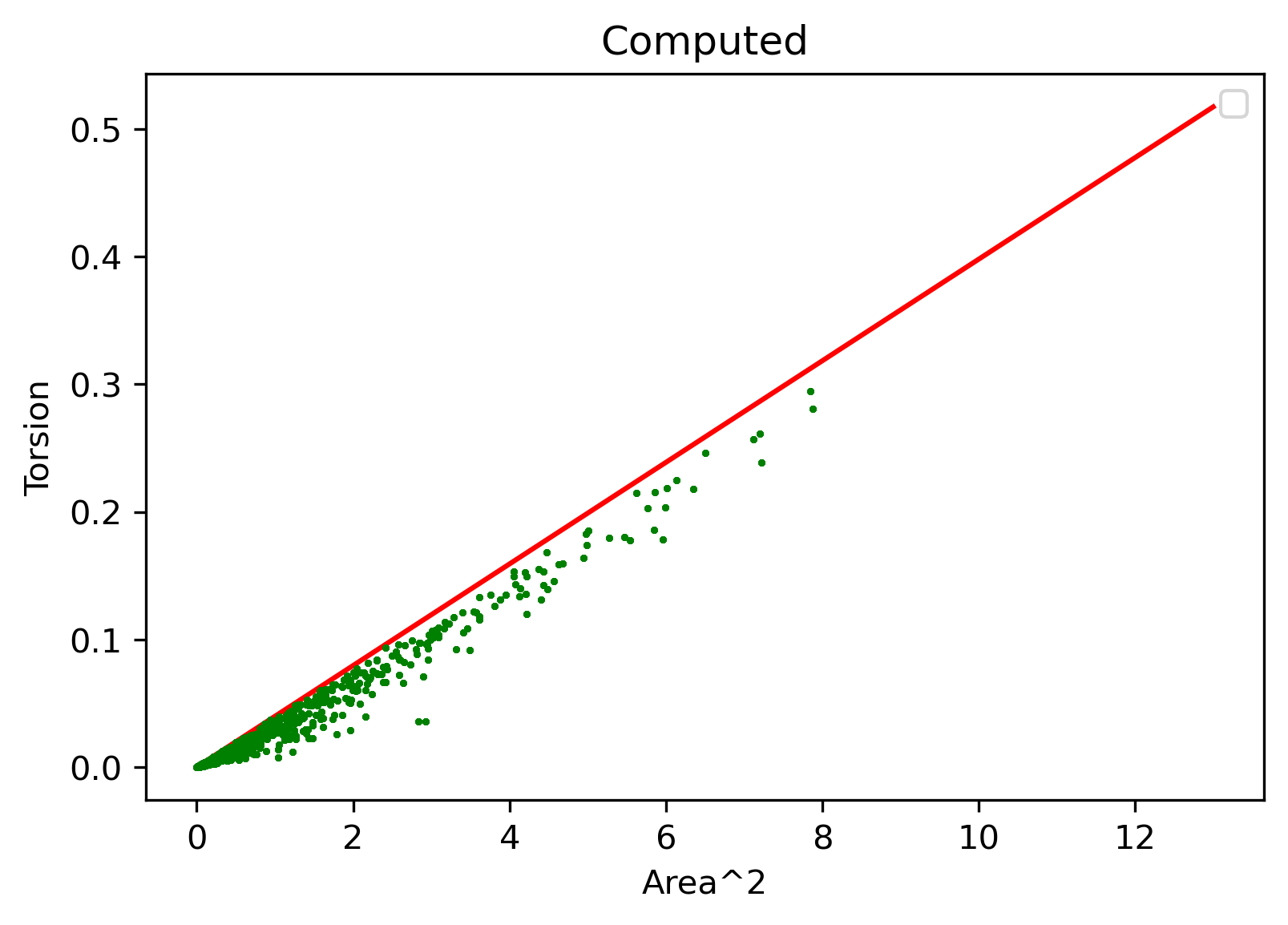}
	\includegraphics[width=0.45\textwidth]{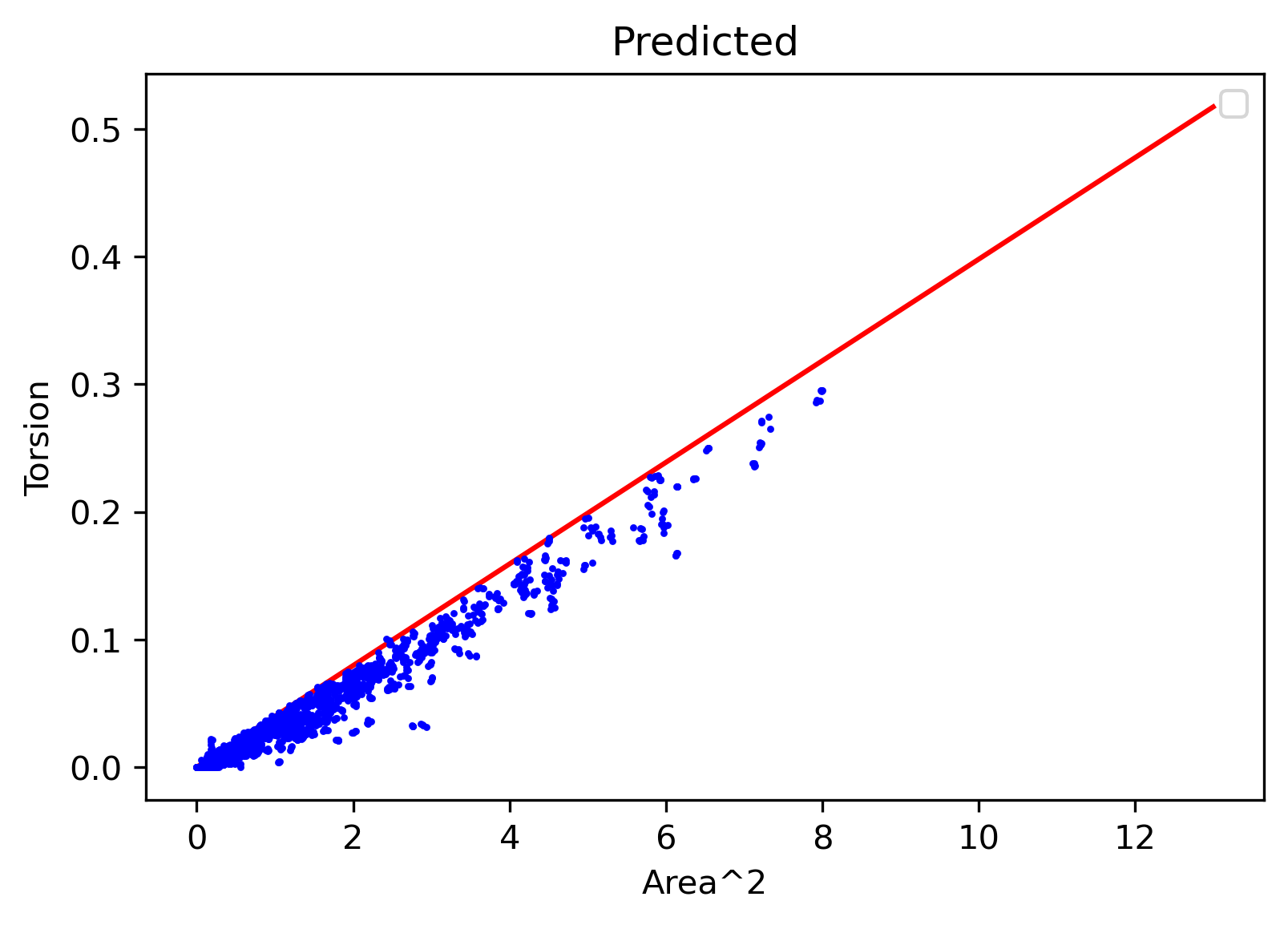}

	\caption{The plot gives the computed value of Torsion (dots) w.r.t the square of the area of the domain. In this plot, test set $\mathrm{Te}$ is considered. The red line corresponds to disks. According to Saint-Venant property all point should lay below such a line.}
	\label{fig:test comparison}
\end{figure}

The values of the loss function and MSE function for the test set are in perfect agreement with those of the validation set, more than satifactory in terms of accuracy. In Figure \ref{fig:test comparison} we have a global picture of the computed and the predicted values. The Torsion is plotted aginst the square of the area. The solid line represents the theoretical line of the torsional rigidity of disks. From the Saint-Venant property \eqref{eq:upperboundBall} no point can lay above that line. The picture clearly shows that the predicted values meet such rule. Apart in few cases for very small values of the Torsion where, the choice of the loss function (mean square error) account for less accuracy.\\

As we already mention, in the loss function we use the standard MSE measure. We give more importance, on purpose, to those set having a larger value of torsion, and we pay the price to lose accuracy on smaller values. But small value of the reference set $Y$ are intrinsecally less accurate, since they mostly correspond to small domains with a less accurate resolution of the boundary.\\

By construction, every random set is enclosed in the box $\mathcal{Q}$, and no value of the Torsion happens to exceed $1$. In order to provide a more precise estimate of the accuracy, we splitted each one of the Training, Validation and Test set, in three subgroups according to the value of the Torsion.
\begin{itemize}
    \item[LV] Large value of the Torsion: $T>0.1$
    \item[SV] Small values of the Torsion $0.01<T<0.1$
    \item[NV] Negligeable values of the Torsion $T<0.01$
\end{itemize}
For LV and SV we computed the Mean Absolute Percentage Error (MAPE):
$$
\frac{1}{M}\sum_{i=1}^M \frac{|y_i-\hat y_i|}{y_i}.
$$
The results are in Table \ref{tab_mape}. The MAPE, as expected, increases on small values of Torsion. But the overall results show that the NN succesfully approximate the functional.

\begin{table}[h]
    \centering
    \begin{tabular}{|l|c|c|} \cline{2-3} 
        \multicolumn{1}{c|}{} & MAPE LV & MAPE SV \\ \hline 
        Tr (Training set) & $3.2\%$ & $6.7\%$   \\ \hline
        Val (Validation set) & $2.8\%$ & $6.7\%$  \\ \hline
        Te (Test set) & $3\%$ & $6.9\%$  \\ \hline
    \end{tabular}
    \caption{MAPE on VGG-16 results. LV stands for Large Values ($T>0.1$). SV stands for Small Values ($0.01<T<0.1$).}
    \label{tab_mape}
\end{table}

\subsection{Polygonal domains}\label{sect:TestPoly}
\begin{figure}
	\centering
	\includegraphics[width=0.8\textwidth]{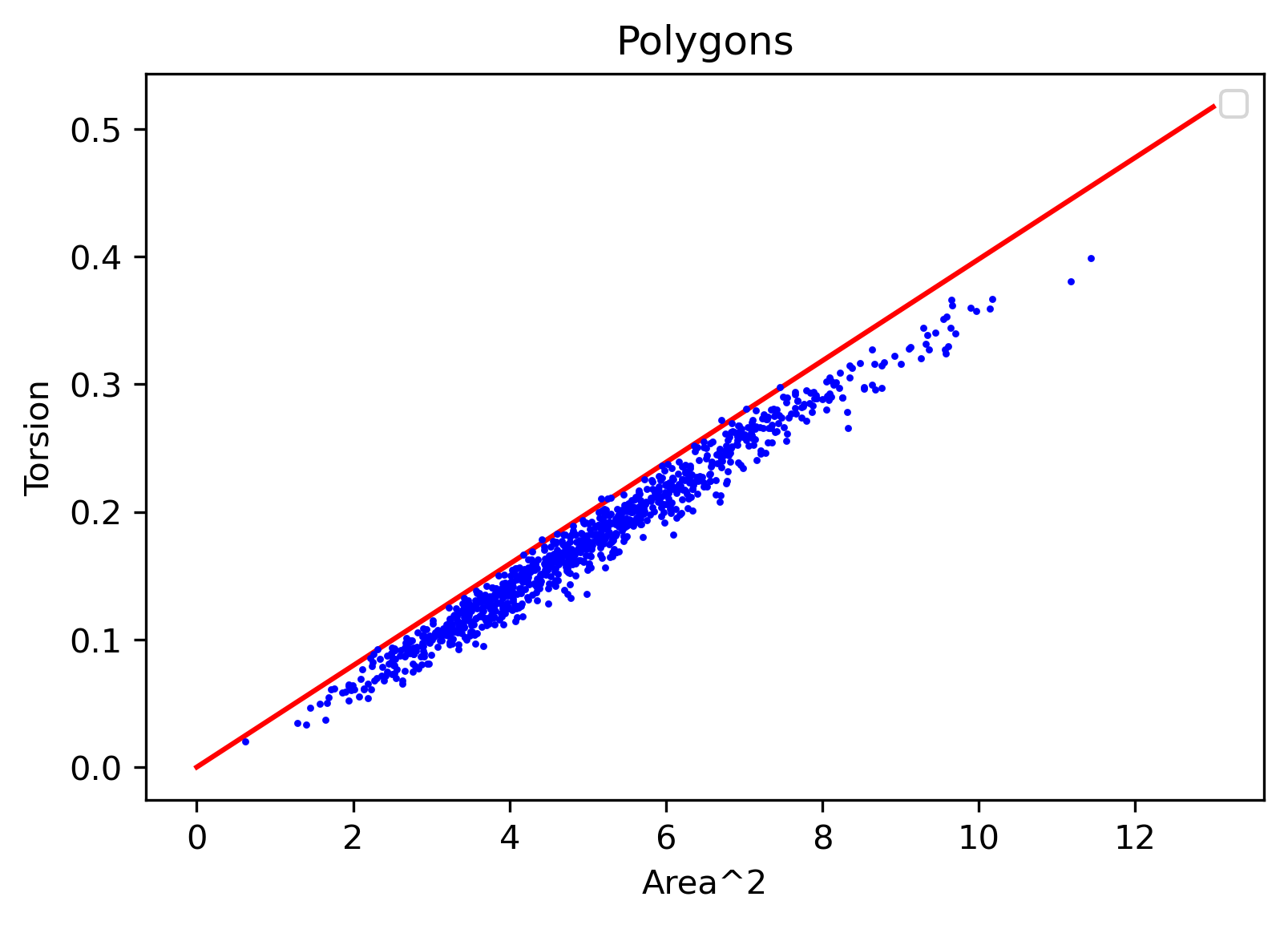}
	\caption{The plot shows predicted values of the Torsion (dots) w.r.t the square of the area of the domains. In this plot, the test set is composed by polygonal domains with number of vertices from $4$ to $11$. In red line is the exact value for disks.}
	\label{fig:poligoni}
\end{figure}

\begin{figure}
	\centering
	\includegraphics[width=0.8\textwidth]{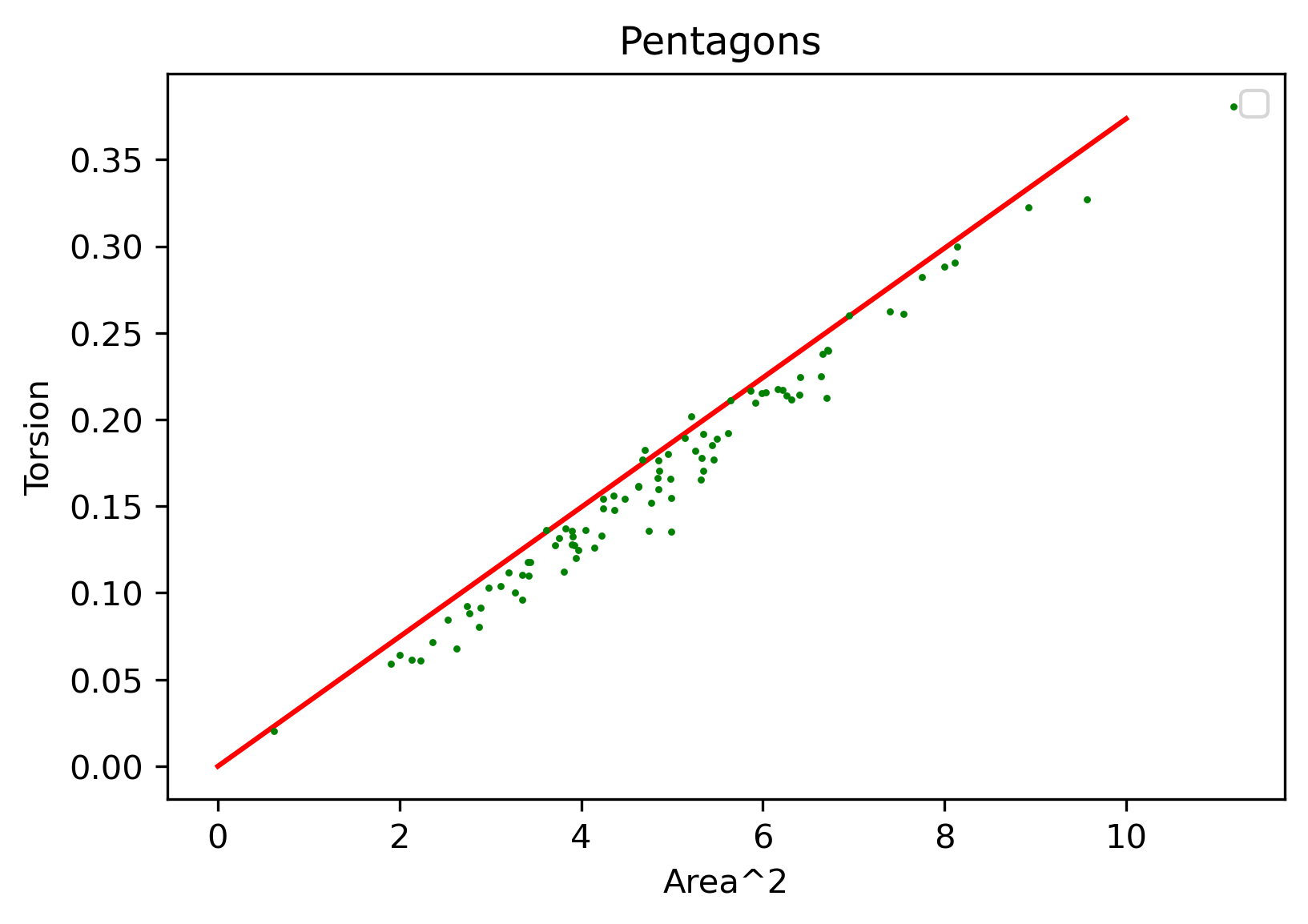}
	\caption{The plot gives the predicted value of Torsion (dots) w.r.t the square of Area of the domain. In this plot, the test set is composed by Pentagons. In red line is the reference numerical value obtained for the regular pentagon.}
	\label{fig:poligoni no regolari 5}
\end{figure}

\begin{figure}
	\centering
	\includegraphics[width=0.8\textwidth]{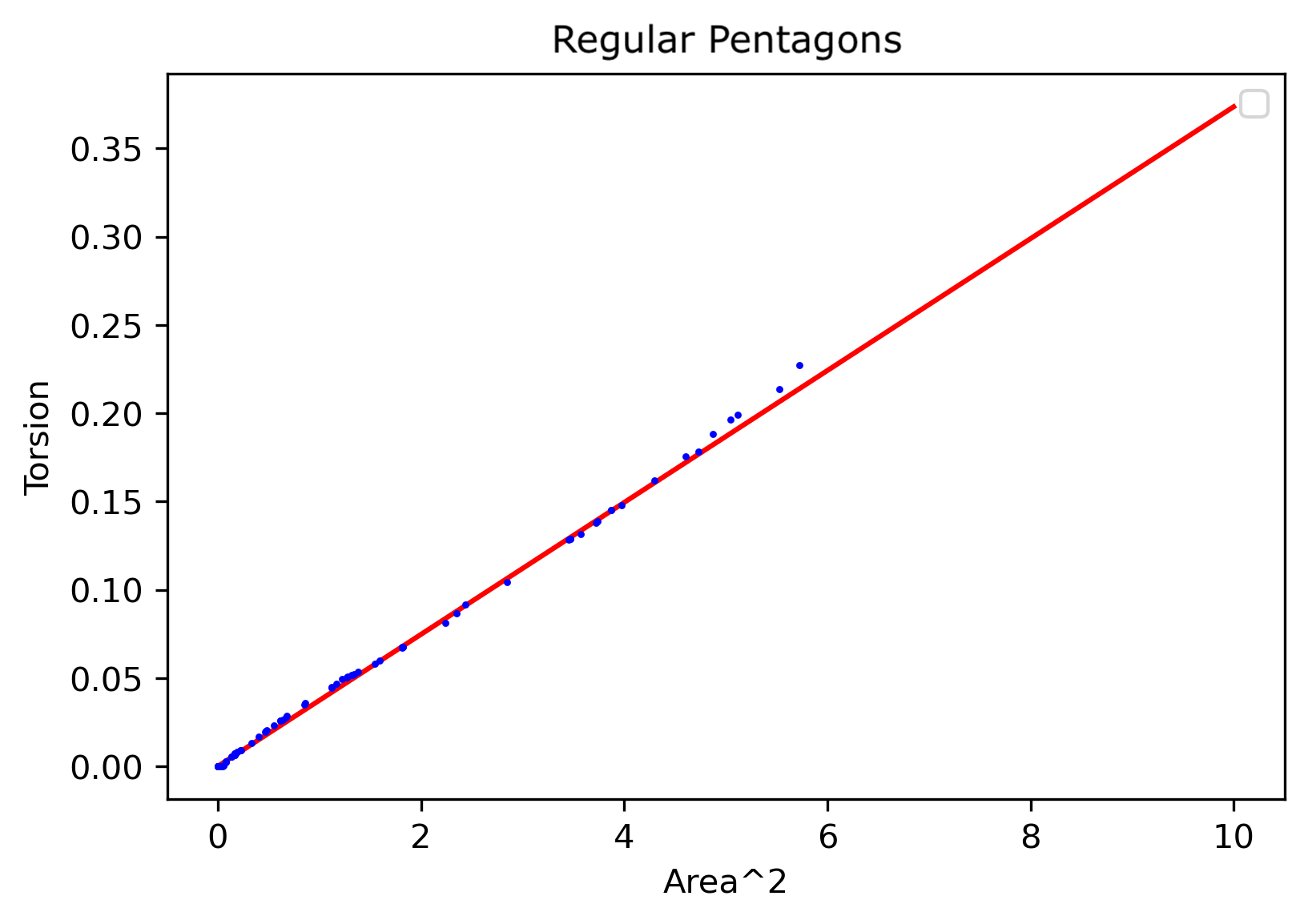}
	\caption{The predicted Torsion on regular pentagons. The slope of the red line comes from numerical computation.}
	\label{fig:pentagoni_r}
\end{figure}

We decided to perform an additional test, restricting the analysis on polygonal domains.
Firstly we generated $1000$ polygonal domains with number of vertices from $4$ to $11$. In Figure \ref{fig:poligoni} we provide a picture of the distribution of the predicted values and we check the agreement with the Saint-Venant property).\\ 
Thereafter we restricted to pentagons. Results are presented in Figure \ref{fig:poligoni no regolari 5}, as usual Torsion is plotted against squared area. According to Conjecture \ref{congettura}, all points in the picture should lay below the continuos line representing regular pentagons. The agreement is almost perfect. Please observe that Conjecture \ref{congettura} is, to a certain extent, a refinement of the Saint-Venant property and therefore requires a greater accuracy.
For the sake of completeness predicted values of the Torsion on regular pentagons are presented in Figure \ref{fig:pentagoni_r}.

\subsection{The case of dilatations}
\begin{figure}
	\centering
	\includegraphics[trim=0 0 0 346, clip, width=0.8\textwidth]{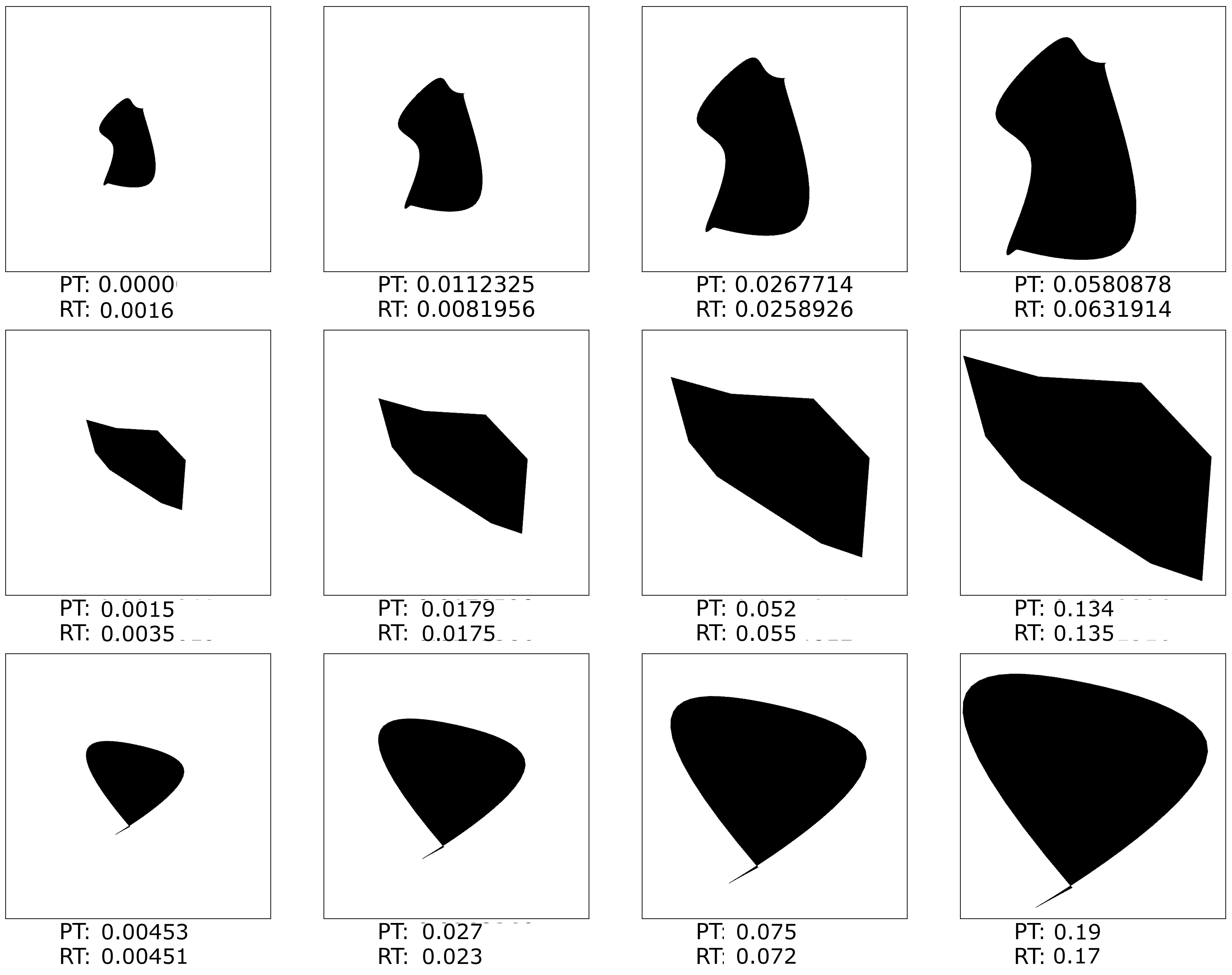}
	\caption{Predicted vs Reference values of the Torsion under dilations.}
	\label{fig:Test_ingrandimento}
\end{figure}

The scaling law under dilation provides yet another way to check a trained network. In Figure \ref{fig:Test_ingrandimento} we present two typical examples. In the first one (top row), the accuracy of the network increases as expected when the set (and the Torsion as well) gets bigger. In the second example (bottom row) there is a very good accuracy on a small scale, more or less preserved under dilation. One can take advantage of the scaling law to look for the estimate of the Torsion on the {\it right} scale. On a small scale ($T<0.01$) one takes the risk of a slack prediction.

\section{Untrained cases}\label{sect:Predictions}
By untrained cases we refer to classes of sets that the DNN has {\it never seen} before. We left this classes unexplored on purpose, to check the extrapolation skills of the DNN. Bear in mind that no matter how large is the training set, there will be always some blind spot in which the DNN has to extrapolate. \\
We consider:
\begin{itemize}
    \item[(A)] Disconnected sets.
    \item[(B)] Multiple connected sets.
\end{itemize}

\subsection{(A) Disconnected sets}
\begin{figure}
	\centering
	\includegraphics[width=0.8\textwidth]{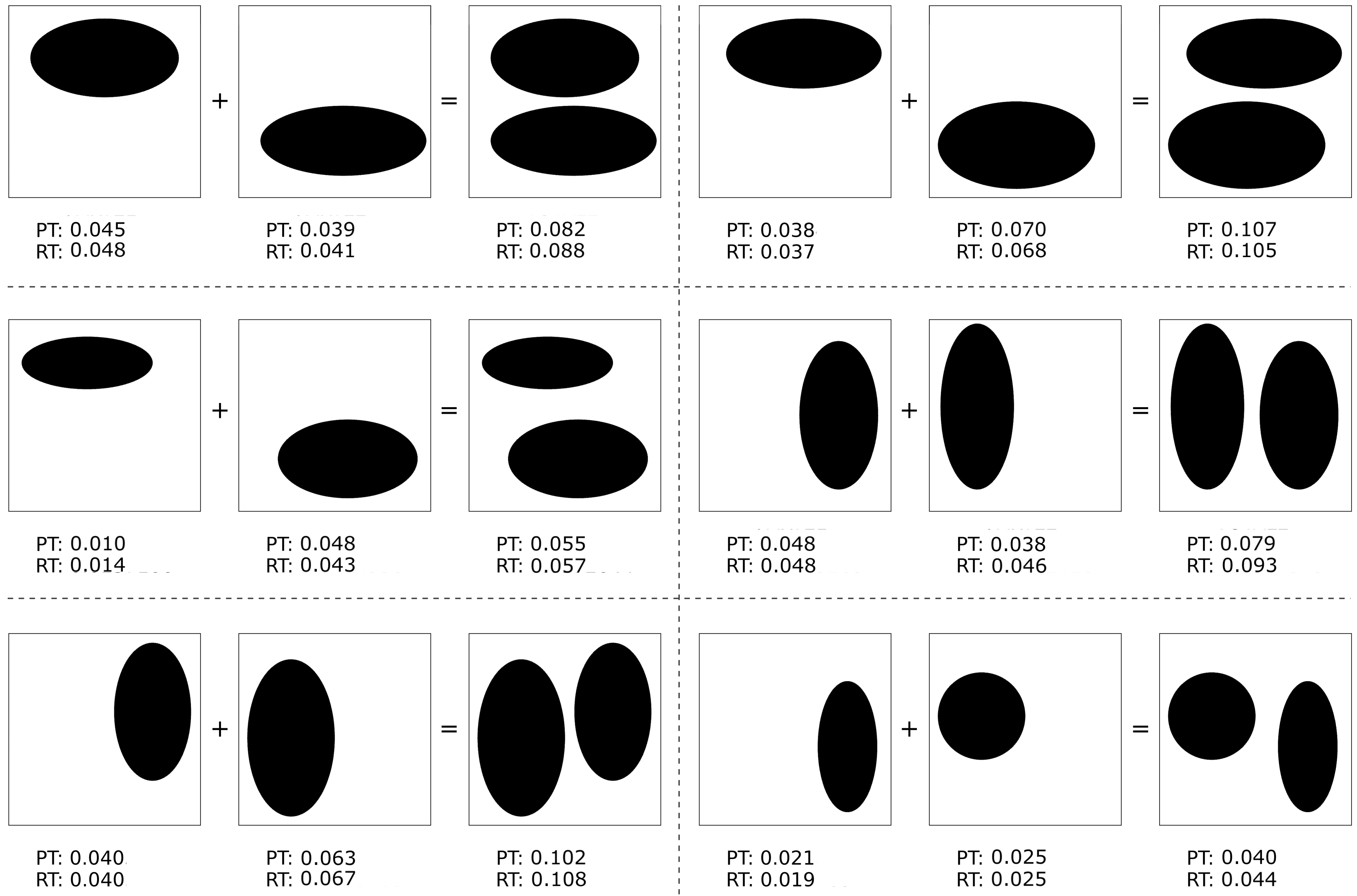}
	\caption{Predicted vs. reference values in the case of union of domains.}
	\label{fig:Union}
\end{figure}

From Section \ref{sect:Problem} we know that on disconnected sets (A) the Torsion is an additive functional: the Torsion of the union is the sum of the torsions. In Figure \ref{fig:Union} we present the results where analytical values and predicted values for unions of ellipses are compared. To the level of accuracy achieved by the DNN we can state that the network {\it understands} that the Torsion is additive under union.\\

\subsection{(B) Multiple connected sets}
\begin{figure}
	\centering
	\includegraphics[trim= 0 300 200 0, clip, width=0.8\textwidth]{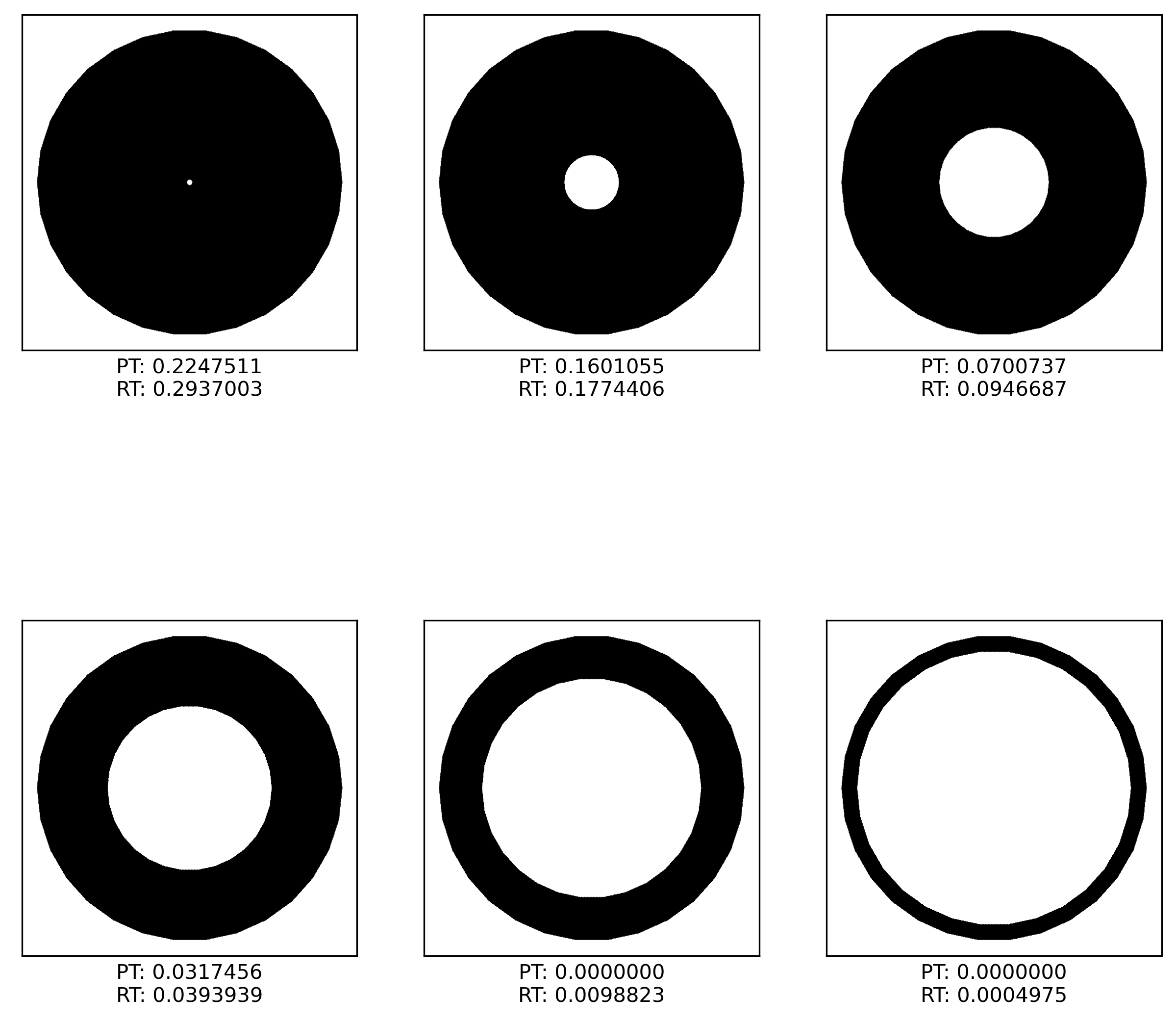}
	\caption{Prediction vs. reference values of the Torsion for two annuli. On the left, despite the small size of the hole, the network ``understands" that the Torsion is noticeably smaller than on the unit disk (where the value is $0.39$). On the right, we show a larger hole, and the network shows greater accuracy.}
	\label{fig:concentrici}
\end{figure}

\begin{figure}
	\centering
	\includegraphics[width=\textwidth]{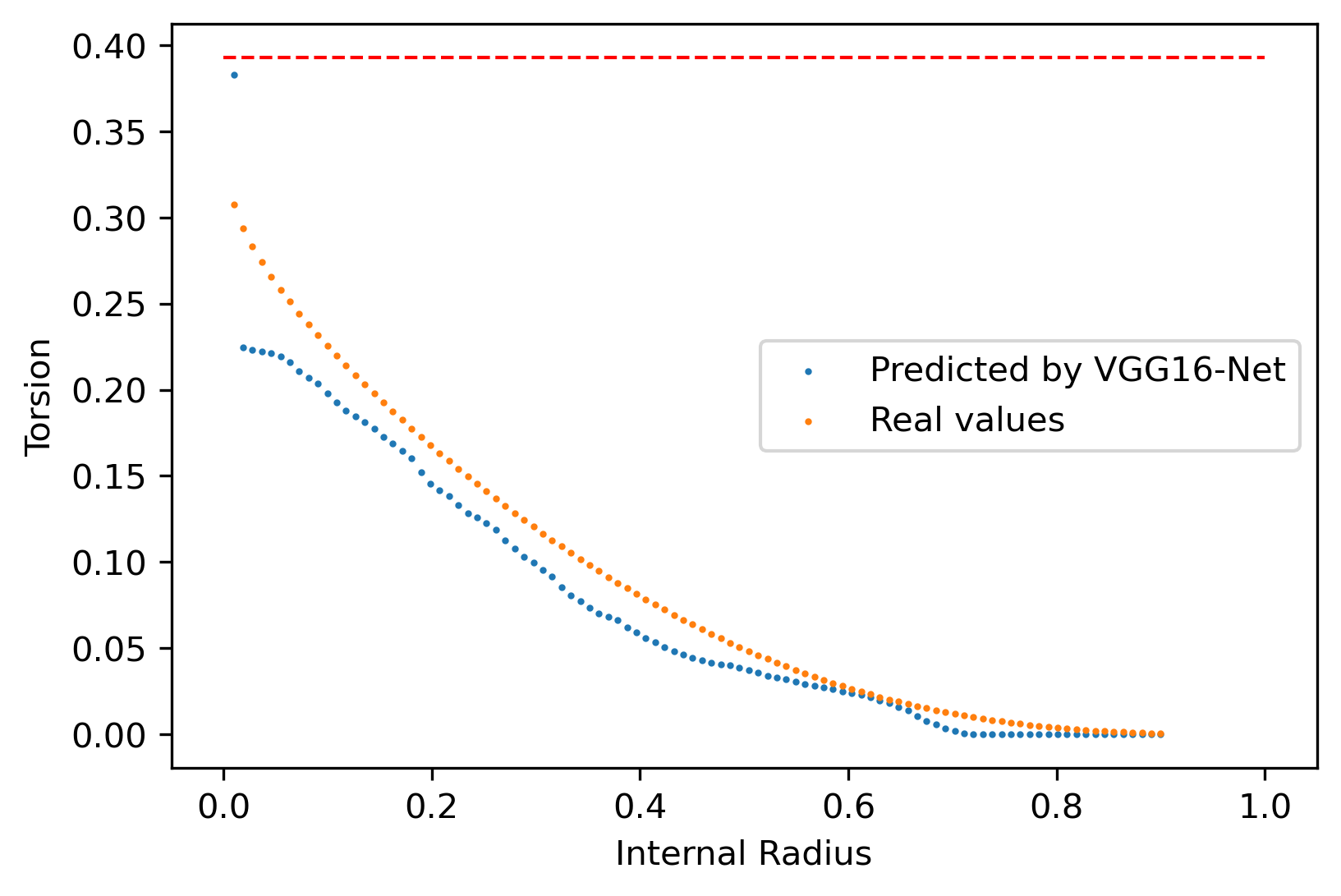}
	\caption{Predicted values of Torsion (blu dots) vs. Reference value of the Torsion (orange dots) for annular rings w.r.t. the internal radius. The dashed red line is the Torsion value of the unitary ball, i.e.: ($\frac{\pi}{8}$)}.
	\label{fig:concentrici2}
\end{figure}

More delicate is the case of multiple connectd sets. The Torsion can be very sensitive to {\it holes}. 
Let us consider the explicit computation of the Torsion on the annular ring $\mathcal{A}_{r1}$ in \eqref{eq:TorAnnularRing} and compare with the Torsion on the disk $\mathcal{D}_1$ \eqref{eq:TorBall}. Following the intuition, $T(\mathcal{A}_{r1})$ converges to $T(\mathcal{D}_1)$ as $r$ goes to $0$. However such a convergence is very slow (logarithmic w.r.t. $r$). If we consider for instance the Torsion of the disk of radius one $T(\mathcal{D}_1)=0.39$ and the Torsion of the annular ring $T(\mathcal{A}_{0.02\ 1})=0.29$, we see that a tiny hole of radius $0.02$ accounts for a decrease of $25\%$. Of course, expecting an accurate prediction with such a small hole is unreasonable, for several reasons, not the least the resolution of the images. It is quite surprising that the DNN catches the behavior (see Figure \ref{fig:concentrici}(a)). As the hole increases the absolute error decreases (see Figure \ref{fig:concentrici}(b)). We provide the general picture in Figure \ref{fig:concentrici2}.
The Network perfectly catches the behavior w.r.t. the internal radius, even if in general it underestimates the value of the Torsion. 

\begin{figure}
	\centering
	\includegraphics[width=\textwidth]{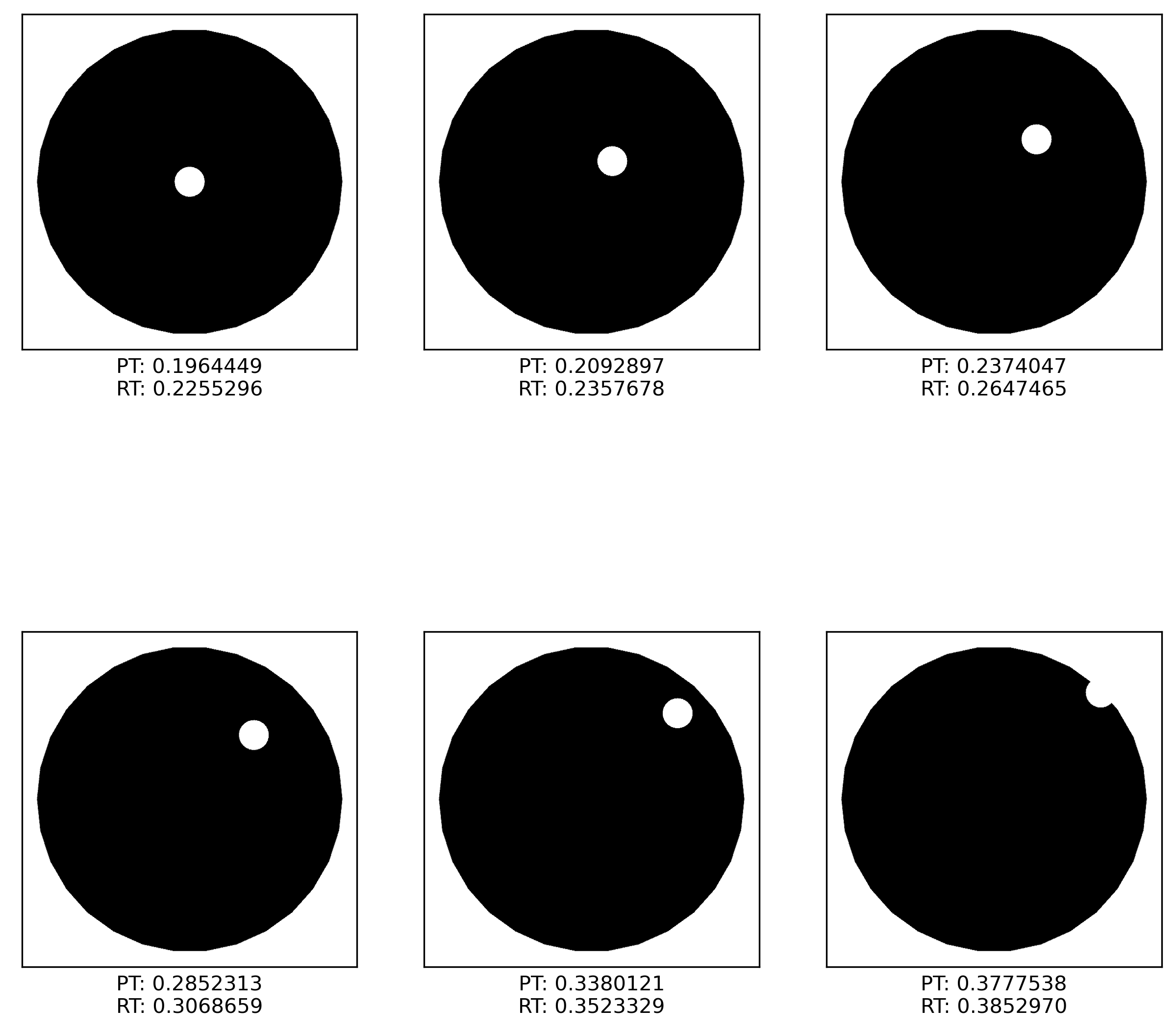}
	\caption{Eccentric annuli. The Torsion increases as the distance between centers increases. Bear in mind that none of the training sets barely resembles such examples.}
	\label{fig:cerchio mobile}
\end{figure}

\begin{figure}
	\centering
	\includegraphics[width=\textwidth]{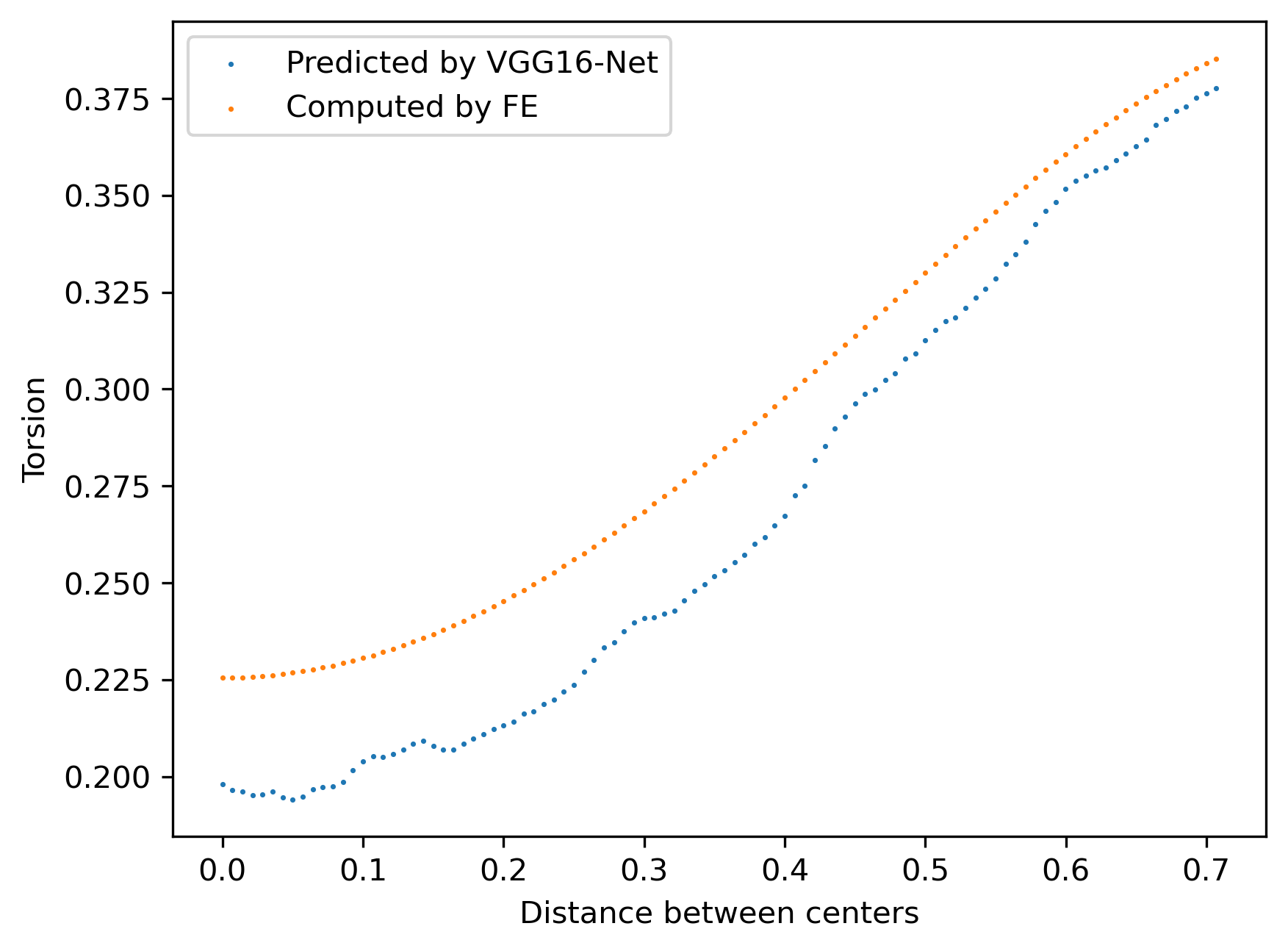}
	\caption{The predicted values of the Torsion (blue dots) w.r.t the distance between centers for eccentric annular rings. The corresponding reference values (numerical computation) are plotted in red.}
	\label{fig:cerchio mobile_2}
\end{figure}

A more intriguing test was performed on eccentric annuli. It is well known that for given inner and outer radii, the more eccentric is the annulus the bigger is the torsional rigidity, the concentric annulus providing the smallest possible value. In Figure \ref{fig:cerchio mobile} we show some examples of the prediction vs. reference.
In Figure \ref{fig:cerchio mobile_2} we provide the plot of the Torsion w.r.t to the distance between the centers (of inner circle and outer circle). With few exceptions the DNN captures the monotonicity. There is a general underestimation of the torsion, but again the network has never seen anything similar during the training.

\section{Conclusions}\label{conclusioni}
We consider this paper as a tentative new approach to the branch of Calculus of Variation known as Shape Optimization. In this paper we have considered the training of a Deep Neural Network for the prediction of Torsional Rigidity.
The use of the Torsion is purely accidental, in the sense that it is nothing but a nice example of Shape Functional. From torsional rigidity to thermal insulation, electrostatic capacity etc., the steps are small, which means that the same approach can be easily generalized in many other meaningful physical problems.
As ususal in Machine Learning, we have trained a NN and checked it on training, validation, and test sets. We noticed a very good accuracy. However, our training data set did not include some classes of domains on purpose. 
Indeed, we checked that the NN properly extrapolate informations also on sets topologically different from those we fed the network during the training.\\
Apart from that, a well trained NN has some undeniable advantage:
\begin{itemize}
    \item Short computational time (realtime answer).
    \item It works with images, even hand drawings, and does not require any analytical expression or approximation of the the boundary of the set.
    \item It can be easily employed in an automated optimization procedure. For instance to minimize (maximize) the shape functional under some geometrical constraints.
\end{itemize} 

Our future perspectives include

\begin{itemize}
    \item Try out new functionals and check other networks architectures to find out the most competitive (in terms of accuracy).
    \item 
    Look for different {\it loss functions} which possibly allow to recover good approximation also at a ``small'' scale.
    \item Train a NN on some Shape functional and then set up a shape optimization problem to be solved by means of NNs, perhaps in conjunction with genetic algorithms.
\end{itemize}

\bibliographystyle{plain}
\bibliography{Biblio.bib}
\end{document}